\documentstyle[12pt,amscd]{amsart}
\newtheorem{prop}{Proposition}
\newtheorem{thm}{Theorem}
\newtheorem{cor}{Corollary}
\newtheorem{lem}{Lemma}

\theoremstyle{definition}

\newtheorem{df}{Definition}

\theoremstyle{remark}
\newcommand{\N}{{\Bbb N}}

\newcommand{\Z}{{\Bbb Z}}
\newcommand{\C}{{\Bbb C}}

\newcommand{\dint}{\displaystyle\int}

\makeatletter
\def\@currentlabel{2.1}\label{e:dispaa}
\def\@currentlabel{2.21}\label{e:dispau}
\def\@currentlabel{2.22}\label{e:dispav}
\def\@currentlabel{2.23}\label{e:dispaw}
\def\@currentlabel{2.24}\label{e:dispax}
\def\theequation{\thesection.\@arabic\c@equation}

\makeatother

\makeatletter
\def\alphenumi{%
  \def\theenumi{\alph{enumi}}%
  \def\p@enumi{\theenumi}%
  \def\labelenumi{(\@alph\c@enumi)}}
\makeatother

\begin{document}

\title{SOME PROPERTIES OF THE SPACE OF COMPACT OPERATORS}
\author{Narcisse Randrianantoanina}
\address{Department of Mathematics, Bowling Green State University,
 Bowling Green, OH 43403}
\email{nrandri@@andy.bgsu.edu}
\maketitle

%\vskip 2in

\begin{abstract}  Let $X$ be a separable Banach space, $Y$ be a Banach
space and $\Lambda$ be a subset of the dual group of a given compact
metrizable abelian group. We prove that if $X^*$ and $Y$ have the type
I-$\Lambda$-RNP (resp. type II-$\Lambda$-RNP) then $K(X,Y)$ has the type
I-$\Lambda$-RNP (resp. type II-$\Lambda$-RNP) provided $L(X,Y)=K(X,Y)$.
Some corollaries are then presented as well as results conserning the
separability assumption on $X$. Similar results for the NearRNP and the
WeakRNP are also presented.
\end{abstract}

\section{INTRODUCTION}
Let $X$ and $Y$ be two Banach spaces. We denote by $L(X,Y)$ (resp.
$K(X,Y)$) the Banach space of all bounded  (resp.
compact bounded) linear operators from
$X$ into $Y$.

This note is devoted to study different types of Radon-Nikodym Property
(RNP) of the space $K(X,Y)$. Recall that in \cite{DM} Diestel and
Morrison proved the following theorem:

\begin{thm}
Suppose that  $X$ is a separable Banach space such that $X^*$ has  the
RNP and $Y$ is a Banach space with the RNP. Suppose in addition that
$L(X,Y) =K(X,Y)$. Then the space $K(X,Y)$ has the RNP.
\end{thm}
Later Andrews \cite{A2} showed that one can remove the separability
assumption on $X$ if either $Y$ is a dual space or if $X$ satisfies the
following topological property: $(^*)$  The weak$^*$-closure of every
bounded norm separable subset of $X^*$ is weak$^*$-metrizable.

In is natural to ask if the same type of result holds for different types
of Radon-Nikodym properties such as Analytic Radon-Nikodym property
(ARNP), Weak Radon-Nikodym property (WeakRNP), Near Radon-Nikodym
property (NearRNP)
(a weakening of the RNP introduced by
Kaufman, Petrakis, Riddle and Uhl \cite{KPRU})
,...

 In the first part of this note we will present some conditions on
$X$ and $Y$ so that the separability of the space $X$ is no longer
needed. Then we consider some types of Radon-Nikodym properties
associated with subsets of countable  discrete abelian group (type
I-$\Lambda$-RNP and type II-$\Lambda$-RNP) which generalize the usual
Radon-Nikodym property and the Analytic Radon-Nikodym property. These
properties were introduced by  Dowling \cite{DO2} and Edgar \cite{EDG}.
We show that the theorem  above still holds if one replace the
usual RNP by the type I-$\Lambda$-RNP or the type II-$\Lambda$-RNP. In
particular Theorem 1. is valid for the Analytic Radon-Nikodym property.
We prove also that similar result can be obtained for the
NearRNP
 and finally we discuss
when a $K(X,Y)$-valued measure has a Pettis-integrable density.
We will show that the above theorem holds for the weakRNP if the range
space $Y$ is a dual space.

All unexplained terminologies can be found in \cite{D1} and \cite{DU}.

\section{THE RADON-NIKODYM PROPERTY FOR THE SPACE OF OPERATORS}
In this section we will provide some sufficient conditions on the Banach
spaces $X$ and $Y$ so that the above theorem is still valid for $X$ non
separable.

We say that a series $\sum\limits_{n=1}^\infty x_n $ is a {\bf weakly
unconditionally
Cauchy (w.u.c)} in $X$ if it satisfies one of the following
equivalent statements:
\begin{enumerate}
 \item [(a)] $\sum\limits_{n=1}^\infty |x^*(x_n)| < \infty$, for every
$x^* \in X^*$;

 \item [(b)] $\sup \{ || \sum\limits_{n \in \sigma} x_n ||; \ \sigma
\ \text{finite subset of}\ {\Bbb N} \} < \infty$;

\item [(c)] $\sup\limits_{n} \sup\limits_{\epsilon_i = \pm 1} ||
\sum\limits_{i=1}^n \epsilon_i x_i || <\infty.$
\end{enumerate}

Some more equivalent formulations can be found in \cite{D1}.

Let us begin by defining the following property introduced by
Pe\l czynski in \cite{PL1}.

\begin{df}(Pe\l czynski) A Banach space $X$ has {\bf property (u)}
if for any weakly Cauchy
sequence $(e_n)_n$ in $X$, there exists a weakly unconditionally  Cauchy
series $\sum\limits_{n} x_n $ in $X$ such that the sequence $(e_n -
\sum\limits_{i=1}^n x_i)$ converges weakly to zero in $X$.
\end{df}

\begin{thm} Let $X$ and $Y$ be Banach spaces such that $X^*$ and $Y$ have
the RNP   and $L(X,Y)=K(X,Y)$. If either $X$ has property (u) or $Y$ is
weakly sequentially complete then $K(X,Y)$ has the RNP.
\end{thm}

The main ingredient for the proof is the following proposition due to
Heinrich and Mankiewicz ( see Proposition 3.4 of \cite{HM}).

\begin{prop}(Heinrich and Mankiewicz)\cite{HM} \label{prop HM}
 Let $X$ be a Banach space and $X_0$ be a separable
subspace of $X$. Then there exists a separable subspace $X_1$ of $X$
containing $X_0$ and an isometric embedding $J: X_{1}^* \to X^*$ with
the property that $\langle z, Jz^* \rangle = \langle z, z^* \rangle $
for every $z \in X_1$ and $z^* \in X_{1}^*$. In particular, $J(X_{1}^*)$
is
1-complemented in $X^*$.
\end{prop}

\begin{prop}\label{thm B}
Let $Z$ be a separable subspace of $X^*$, then there exist  a
1-complemented subspace $Z_1$ of $X^*$ with $Z \subset Z_1
\subset X^*$ and a separable subspace $X_1$ of $X$ suth that $Z_1$ is
isometric to $X_{1}^*$.
\end{prop}

\begin{pf} If $\{f_n , n \geq 1\}$ is
a countable
dense subset of the unit ball of $Z$  and $\{x_{n,j};\  n,j \geq 1 \}$
is a sequence in $X$ such that for every $n \in \N$,  $\lim\limits_{j \to
\infty } f_n (x_{n,j}) =||f_n||$.
 Let $X_0$ be the separable subspace of $X$ spanned by the sequence
$(x_{n,j})_{n,j}$. By Proposition \ref{prop HM}, there exist  a
separable subspace
$X_1$ of
$X$
with  $X_0 \subset X_1 \subset X$ and $J: X_{1}^* \to X^*$ as in
Proposition
\ref{prop HM}. We will show that $Z$ is isometrically isomorphic
to  the subspace
$J(X_{1}^*)$  of $X^*$.
 For that let
$Q: X^* \to X_{1}^*$ be the restriction map and $i: Z \to X^*$ the
inclusion . We claim that $Q\circ i: Z \to X_{1}^*$ is an isometry. In
fact since $X_0$ is norming for $Z$, we have for every $f \in Z$,
  $$ ||Q\circ i (f)||= \sup_{||z||\leq 1} \langle Q(f),z \rangle
=\sup_{||z|| \leq 1} \langle f,z \rangle = ||f|| $$
which shows that $Q\circ i$ is an isometry. Now since  $J$ is an
isometry,
we have $Z$ embedded isometrically into $J(X_{1}^*)$ and the proposition
is proved.
\end{pf}

We are now ready to provide the proof of the theorem: It is enough to
show  that every separable subspace $\cal S$ of $K(X,Y)$ is isometric to
a subspace of $K(X_1,Y_1)$ where $X_1$ is a separable subspace of $X$ and
$Y_1$ is a separable subspace of $Y$ and such that $L(X_1,Y_1)
=K(X_1,Y_1)$.

Let $\cal S$ be a separable subspace of $K(X,Y)$. It is clear by
compactness that  the space $\left\{Tx ;\ T\in \cal{S},\ x \in X
\right\}$ is
separable.  The space $Y_1 = \overline{\text{span}}\left\{Tx ;\ T\in
\cal{S},\
x \in X \right\}$ is separable and $\cal{S} \subset K(X,Y_1)$. Using
similar argument with the adjoints, we get that
 the space $Z= \overline{\text{span}}\left\{T^*y^* ;\ T\in
\cal{S},\   y^*  \in Y^* \right\}$ is a separable subspace of $X^*$.
Let $Z_1$ and $X_1$ as in Proposition 2. It is clear that the restriction
map from $\cal{S}$ into $K(X_1,Y_1)$ is an isometry and we claim that
$L(X_1,Y_1)=K(X_1,Y_1)$. For that let $J: X_{1}^* \longrightarrow X^*$
and $Q: X^* \longrightarrow X_{1}^*$ as before and fix $\theta \in
L(X_1,Y_1)$. Assume that $X$ has property (u); since $X_1$ does not
contain
$\ell^1$ and $Y_{1}$ does not contain $c_0$,
by \cite{PL1}, the operator $\theta$ is weakly compact. Similarly if
$Y_1$ is weakly sequentially complete, the operator $\theta$ is weakly
compact. So in both cases, the operator $\theta$ is weakly compact.
Consider $J\circ \theta^*: Y_{1}^* \longrightarrow X^*$; by
the weak compactness of $\theta$, we have
  $$(J\circ \theta^*)^* (X^{**})= \theta^{**}\circ J^*(X^{**}) \subset
Y_1$$
and since $L(X,Y)=K(X,Y)$, we get that $(J \circ\theta^*)^* |_{X}$ is
compact and therefore $(J\circ\theta^*)^*(B_{X})$ is relatively compact
in $Y_1$.
     $$(J\circ\theta^*)^*(B_{X^{**}}) \subset
\overline{(J\circ\theta^*)^*(B_{X})}^{||\  ||} $$
which shows that $ J\circ\theta^*)^*(B_{X^{**}})$ is relatively compact.
Hence $J\circ\theta^*$ is compact. To complete the proof, let $\pi: X^*
\longrightarrow J(X_{1}^*)$ be the norm 1- projection; it can be easily
checked that
$\theta^* =J^{-1}\circ\pi\circ(J\circ\theta^*)$ which proves that
$\theta^*$ (and hence $\theta$) is compact. The theorem is proved.

\section{SOME VARIANTS OF THE RADON-NYKODYM PROPERTY FOR THE SPACE OF
OPERATORS}

Throughout the remaining of this paper
$G$ will denote a compact metrizable abelian group, $\cal B(G)$ is the
$\sigma$-algebra of the Borel subsets of $G$, and $\lambda$ the
normalized Haar measure on $G$. We will denote by $\Gamma$ the dual group
of $G$ i.e the set of continuous homomorphisms $\gamma:G \longrightarrow
\C$  ($\Gamma$ is a countable discrete abelian group).

\noindent Let $X$ be a Banach space and $1 \leq p \leq \infty$, we will
denote by $L^p(G,X)$ the usual Bochner function spaces for the measure
space $(G,\cal B (G), \lambda)$, $M(G,X)$
the space of $X$-valued countably additive measure of bounded variation,
$C(G,X)$ the space of $X$-valued continuous functions and
$M^\infty(G,X)= \{ \mu \in M(G,X) , |\mu| \leq C\lambda \ \text{for some
$C >0$} \}$.
\begin{itemize}
\item [(i)] If $f \in L^1(G,X)$, we denote by $\hat f $ the Fourier
transform of $f$ which is the map from $\Gamma$ to $X$ defined by $\hat f
(\gamma)= \int\limits_{G} \bar\gamma f \ d\lambda$.

\item [(ii)] If $\mu \in M(G,X)$, we denote by $\hat\mu$ the Fourier
transform of $\mu$ which is the map from $\Gamma$ to $X$ defined by
$\hat\mu (\gamma) =\int\limits_{G} \bar\gamma \ d\mu$.
\end{itemize}

\noindent If $\Lambda \subset \Gamma$ is a set of characters, let
\begin{align*}
 L_{\Lambda}^p (G,X) &=\{ f \in L^p (G,X),\ \hat f(\gamma)=0 \ \text{for
all $\gamma \notin \Lambda$}\} \\
 C_\Lambda (G,X) &=\{ f \in C(G,X),\ \hat f(\gamma)=0 \ \text{for all
$\gamma \notin \Lambda$}\} \\
 M_\Lambda (G,X) &=\{ \mu \in M(G,X),\ \hat \mu (\gamma)=0 \ \text{for
all $\gamma \notin \Lambda$}\} \\
 M_{\Lambda}^\infty (G,X) &=\{ \mu \in M^\infty (G,X),\ \hat \mu
(\gamma)=0 \ \text{for all $\gamma \notin \Lambda$}\}
\end{align*}

\begin{df}
 (i) A subset $\Lambda$ of $\Gamma$ is a Riesz set
if and only if $M_\Lambda (G,\C) =L_{\Lambda}^1 (G,\C)$

(ii) A subset $\Lambda$ of $\Gamma$ is a Rosenthal set if and only if
$C_{\Lambda}(G) = L_{\Lambda}^\infty (G)$.
\end{df}

\noindent The following properties where introduced by Edgar \cite{EDG},
and Dowling \cite{DO2}.

\begin{df}
 (i) A Banach space $X$ is said to have type
I-$\Lambda$-Radon Nikodym Property (type I-$\Lambda$-RNP) if and only if
$M_{\Lambda}^\infty(G,X)= L_{\Lambda}^\infty(G,X)$.

(ii) A Banach space $X$ is said to have  type II-$\Lambda$-Radon Nikodym
Property (type II-$\Lambda$-RNP) if and only if
$M_{\Lambda,ac}(G,X)=L_{\Lambda}^1 (G,X)$ where $$M_{\Lambda,ac}(G,X) =\{
\mu \in M_{\Lambda}(G,X), \ \text{$\mu$ is absolutely continuous with
respect to $\lambda$}\}.$$
\end{df}

\noindent {\bf Remarks:} (a) It is obvious that type II-$\Lambda$-RNP
implies type I-$\Lambda$-RNP.

(b) Since $\cal B (G)$ is countably generated, one can see  that  these
two properties are separably determined.

(c) If $G={\Bbb T}$ then $\Gamma =\Z$. Then type I-$\Z$-RNP is
equivalent
to type II-$\Z$-RNP which is also equivalent to the usual RNP. Similarly,
type I-$\N$-RNP is equivalent to type II-$\N$-RNP and is equivalent to
the Analytic Radon Nikodym Property (see \cite{EDG}).

(d) If $\Lambda$ is a Riesz subset, then $M_{\Lambda,ac}(G,X) =
M_{\Lambda}(G,X)$.

We are now ready to present our results.
\begin{thm}
Let $X$ and $Y$ be Banach spaces such that:
\begin{itemize}
\item [(i)] $X$ is separable;
\item [(ii)] $X^*$ and $Y$ have type I-$\Lambda$-RNP (resp. type
II-$\Lambda$-RNP);
\item [(iii)] $L(X,Y)=K(X,Y)$;
\end{itemize}
Then $K(X,Y)$ has type I-$\Lambda$-RNP (resp. type II-$\Lambda$-RNP).
\end{thm}

We will present the proof for type I-$\Lambda$-RNP case (the type
II-$\Lambda$-RNP case can be done with minor changes).

Consider $\cal{B}(G)$ the $\sigma$-Algebra generated by the Borel subsets
of $G$ and fix a measure
      $F:\cal B(G) \longrightarrow K(X,Y)$
such that
\begin{itemize}
\item[a)] $|F| \leq \lambda$;
\item [b)] $\widehat{F}(\gamma)=0$\ for \ $\gamma \notin \Lambda$.
\end{itemize}

Our main goal is to show that the measure $F$ has  a Bochner integrable
density.

For $x \in X$, we will denote by $F(.)x$ the $Y$-valued measure
$A \longrightarrow F(A)x$ and similarly, for $y^* \in Y^*$, $F(.)^*y^*$
will be the $X^*$-valued measure $A \longrightarrow F(A)^*y^*$.

Let us begin with the following simple observation:

For every $x \in X$, $y^* \in Y^*$ and $\gamma \in \widehat{G}=\Gamma$,
we have
\begin{align*}
\langle \widehat{F}(\gamma), x \otimes y^* \rangle &= \langle
\widehat{F(.)x}(\gamma),y^* \rangle \\
&= \langle \widehat{F(.)^*y^*}(\gamma),x \rangle
\end{align*}
These equalities imply that for any $x \in X$ and $y^* \in Y^*$,
$F(.)x \in M_{\Lambda}^\infty(G,Y)$  and
$F(.)^*y^* \in M_{\Lambda}^\infty (G,X^*)$.

Notice also that without loss of generality we can and do assume that $Y$
is separable.

We need the following definition for the next Proposition.

\begin{df} Let $E$ and $F$ be Banach spaces and $(G,\cal(B)(G),\lambda)$
be a measure space. A map $T: G\to L(E,F)$ is said to be strongly
measurable if $\omega \to T(\omega)e$ is $\lambda$-measurable for every
$e \in E$.
\end{df}

\begin{prop}\label{lem1}
There exists a strongly mesurable map $\omega \longrightarrow T(\omega)$
$(G \longrightarrow K(X,Y))$ such that:
\begin{itemize}
\item[(a)] $F(A)x= \text{Bochner}-\dint_{A} T(\omega)x\ d\lambda(\omega)$
for every $A \in \cal{B}(G)$ and $x \in  X$;
\item[(b)]  For every $y^* \in Y^*$, the map $\omega \longrightarrow
T(\omega)^*y^*$ is norm-measurable and
$F(A)^*y^* = \text{Bochner}-\dint_{A} T(\omega)^*y^* \ d\lambda(\omega)$
for every $A \in \cal{B}(G)$.
\end{itemize}
\end{prop}

\begin{pf}
Using similar argument as in \cite{DM}, one can construct a strongly
measurable map $\omega \longrightarrow T(\omega)$ such that $F(A)x=
\text{Bochner}-\dint_{A} T(\omega)x \ d\lambda(\omega)$  for every $x \in
X$ and $A \in \cal{B}(G)$ so the first part is proved.

For the second part, notice from the strong measurability of $T(.)$ that
the map $\omega \longrightarrow T(\omega)^*y^*$ is weak$^*$-scalarly
measurable and
 $F(A)^*y^* = \text{weak$^*$}- \dint_{A} T(\omega)^*y^* \
d\lambda(\omega)$ for every $A \in \cal{B}(G)$.

 Since the measure $F(.)^*y^*$ belongs to $M_{\Lambda}^\infty(G,X^*)$ and
$X^*$ has the type I-$\Lambda$-RNP, there exists  a Bochner integrable
map $h_{y^*}: G \longrightarrow X^*$ such that $F(A)^*y^* =
\text{Bochner}-\dint_{A} h_{y^*}(\omega)\ d\lambda(\omega)$ for every $A
\in \cal{B}(G)$ so we get that
$$\int_{A} \langle T(\omega)^*y^*,x \rangle \ d\lambda(\omega)=\int_{A}
\langle h_{y^*}(\omega),x \rangle \ d\lambda(\omega)$$
for every $x \in X$ and $A \in \cal{B}(G)$.
Now fix $\{x_n, n \in \N\}$ a countable dense subset of $X$. There exists
a measurable subset $G^\prime$ of $G$ with $\lambda(G \setminus
G^\prime)=0$ and such that for $\omega \in G^\prime$ and $n \in \N$, we
have  $\langle T(\omega)^*y^*,x_n \rangle = \langle h_{y^*}(\omega),x_n
\rangle$ which of course implies that
 $$T(\omega)^*y^* = h_{y^*}(\omega) \quad \text{for} \ \omega \in
G^\prime$$

so $T(.)^*y^* $ is norm measurable and satisfies the required property.
\end{pf}

\begin{prop} \label{lem2} Let $X$ and $Y$ be separable Banach spaces and
$Z$ be a separable subspace of $X^*$ then the set
$$\cal{A}=\left\{ T \in K(X,Y);\ T^*y^* \in Z, y^* \in Y^* \right\}$$
is separable in $K(X,Y)$.
\end{prop}
\begin{pf}
Let $\Delta$ be the unit ball of $Y^*$ with the weak$^*$-topology. Since
$Y$ is separable $\Delta$ is a compact metric space . Let $C(\Delta)$ be
the Banach space of all continuous functions on $\Delta$ with the usual
sup norm. It is well known that the Banach space $Y$ embedds
isometrically into $C(\Delta)$. Let $J: Y \longrightarrow C(\Delta)$ be
the natural isometry. Consider an operator
$J^\#: K(X,Y) \longrightarrow K(X,C(\Delta))$ defined as follows:
      $$J^\#(T)= J\circ T.$$
It is clear that $J^\#$ is an isometry and  we will show that
$J^\#(\cal{A})$
is separable. For that notice that $J^\#(\cal{A})$ is a subset of
$$\cal{M}=\left\{S \in K(X,C(\Delta));\ S^*u^* \in Z, \ u^* \in
{C(\Delta)}^* \right\}.$$
In fact for every $u^* \in {C(\Delta)}^*$ and $T \in \cal{A}$,  we have
$J^\#(T)^*u^* =(J\circ T)^*u^*= T^*(J^*u^*) \in Z$.
Now since $M(\Delta)={C(\Delta)}^*$ has the approximation property, we
have $$ K(X,C(\Delta))=K_{w^*}(X^{**},C(\Delta))= X^*
\hat{\otimes}_\epsilon C(\Delta)$$
where $K_{w^*}(X^{**},C(\Delta))$ denotes the
space of compact operators from $X^{**}$ into
$C(\Delta)$ that are weak$^*$ to weak
continuous and $\hat{\otimes}_{\epsilon}$ is the injective tensor
product.
Let $S \in \cal{M}$; since $S^* \in K_{w^*}(M(\Delta),Z)= C(\Delta)
\hat{\otimes}_\epsilon Z$, it can be approximated by elements of
$C(\Delta) \otimes Z$ and by duality $S$ can be approximated by elements
of $Z \otimes C(\Delta)$ and therefore $\cal{M} \subset Z
\hat{\otimes}_\epsilon C(\Delta)$ which is a separable space. We are
done
\end{pf}

To complete the proof of the theorem, let us choose a sequence $(A_n)_n
\subset \cal{B}(G)$ such that $\left\{F(A_n),\ n \geq 1, \ A_n \in
\cal{B}(G) \right\}$ is dense in the range of the measure $F$ (this is
possible because $\cal{B}(G)$ is countably generated).

For each $n \geq 1$, the operator $F(A_n)$ is compact so the set
$F(A_n)^*(B_{Y^*})$ is compact in $X^*$ and therefore $F(A_n)^*(Y^*)$ is
separable in $X^*$. Define a subspace $Z$ of $X^*$ as follows:
 $$ Z = \overline{\text{span}}\left\{ \bigcup_{n\geq 1} F(A_n)^*(Y^*)
\right\}.$$
The space $Z$ is obviously a separable subspace of $X^*$ and for every $A
\in \cal{B}(G)$ and $y^* \in Y^*$, we have $F(A)^*y^* \in Z$.
We will need the following lemma.
\begin{lem} \label{lem3}  Let $\cal{A}$ be the separable subspace of
$K(X,Y)$ as in
Proposition \ref{lem2}. For a.e $\omega \in G$, $T(\omega) \in \cal{A}$.
\end{lem}
\begin{pf*}{Proof of Lemma \ref{lem3}}
Let $\{y_{n}^*,\ n \geq 1\}$ be a countable weak$^*$-dense subset of
$B_{Y^*}$. Let $n \in \N$ fixed. Since $F(A)^*y_{n}^* \in Z$ and
$F(A)^*y_{n}^* =\text{Bochner}-\dint_{A} T(\omega)^*y_{n}^* \
d\lambda(\omega)$ for every $A \in \cal{B}(G)$, we get that
$T(\omega)^*y_{n}^* \in Z$ for a.e $\omega \in G$. There exists a
measurable subset $O_n$ of $G$  with $\lambda(O_n)=0$ and for $\omega
\notin O_n$, $T(\omega)^*y_{n}^* \in Z$. Let $O=
\bigcup\limits_{n=1}^\infty O_n$; $\lambda(O)=0$ and if $\omega \notin
O$, we have
   $$T(\omega)^*y_{n}^* \in Z \qquad \forall n\in \N.$$
Now for $y^* \in B_{Y^*}$, choose a sequence $(y_{n_j}^*)_{j \in \N}$
that converges
to $y^*$ for the weak$^*$-topology. Since $T(\omega)^*$ is weak$^*$to
norm continuous, the sequence $\{T(\omega)^*(y_{n_j}^*)\}_{j \in \N}$
converges to $T(\omega)^*y^*$ for the norm-topology  which implies that
$T(\omega)^*y^* \in Z$ for every $\omega \in G \setminus O$ \
(independent
of $y^*$) hence $T(\omega) \in \cal{A}$ for every $\omega \in G \setminus
O$ and the lemma is proved.
\end{pf*}

We complete the proof by noticing that the map $T: G \longrightarrow
K(X,Y)$ that takes
$\omega$ to $T(\omega)$ is $\lambda$-essentially separably valued and
there exists a norming subset of $K(X,Y)^*$ (namely $X {\otimes}Y^*$)
such that the map $\omega \longrightarrow \langle T(\omega),\phi \rangle$
is $\lambda$-measurable for every $\phi \in X \otimes Y^*$ and by the
Pettis-measurability Theorem (see \cite{DU}, Theorem II-2), the map
$\omega \longrightarrow T(\omega)$ is norm-measurable and
it is now clear that the measure $F$ is represented by the map
$T: G \longrightarrow K(X,Y)$.  The proof is complete. \qed

\noindent
{\bf Remark:} The argument used by Andrews in \cite{A2} can be adjusted
to show that for the case where
$Y$ is a dual space, the assumption that
$X$ is separable may be dropped.

 Some corollaries are now in order

\begin{cor} If $X$ is a Banach space with the Schur property and
$\Lambda$ is a Riesz subset then
$L^1(\lambda) \hat{\otimes}_\epsilon X$ (the completion of the space of
Pettis representable measures with the semi-variation norm) has type
I-$\Lambda$-RNP (resp. type II-$\Lambda$-RNP) if and only if $X$ has type
I-$\Lambda$-RNP (resp. type II-$\Lambda$-RNP).
\end{cor}

\begin{pf} Notice that $L^1(\lambda) \hat{\otimes}_{\epsilon} X$ is a
subspace of $M(G) \hat{\otimes}_\epsilon X$ which is isometrically
isomorphic to the space $K(C(G),X)$. Now if $X$ has the Schur property,
$L(C(G),X)=K(C(G),X)$. We are done.
\end{pf}

\begin{cor} Let $X$ be a Banach space having the type I-$\Lambda$-RNP
(resp. type II-$\Lambda$-RNP) and denote by $\ell^1 [X]$ the Banach space
of all W.U.C. series in $X$ normed by
 $$||(x_n)||= \sup\left\{\sum_{n=1}^\infty |x^*(x_n)|;\ x^* \in X^*,\
||x^*|| \leq 1\right\}.$$
Then  $\ell^1[X]$ has the type I-$\Lambda$-RNP (resp. type
II-$\Lambda$-RNP).
\end{cor}
\begin{pf} This is due to the well known fact that $\ell^1[X] =L(c_0,X)$
and since $X$ does not contain any copy of $c_0$, we have
$L(c_0,X)=K(c_0,X)$. An appeal to Theorem 2 completes the proof.
\end{pf}

For the next corollary, let us introduce a new compact
metrizable abelian group $\widetilde G$ which is not necessarily the same
as $G$.
We will denote by $\widetilde\Gamma$ its dual  and $\widetilde\lambda$
its normalized Haar measure. The following result was first obtained by
Dowling \cite{DO3} for the usual RNP. It was also proved in \cite{RS3}
(see also \cite{R}) under the assumption that $\Lambda$ is a Riesz
subset.
\begin{cor} Assume that
\begin{itemize}
\item[(1)] $\widetilde\Lambda$ is a Rosenthal subset of
$\widetilde\Gamma$;
\item[(2)] $X$ is a Banach space with the Schur property and has the type
II-$\Lambda$-RNP.
\end{itemize}
Then the space $C_{\widetilde\Lambda}(\widetilde G,X)$ has the type
II-$\Lambda$-RNP.
\end{cor}

\begin{pf} If $\widetilde\Lambda$ is a Rosenthal subset, then
$C_{\widetilde\Lambda}(\widetilde G)$ has the RNP and
$C_{\widetilde\Lambda}(\widetilde
G)=L_{\widetilde\Lambda}^\infty(\widetilde G) = \left(
L^1(\widetilde G)/L_{{\widetilde\Lambda}^\prime}^1 (\widetilde
G)\right)^*$ is a dual space. Now since $X$ has the Schur property and
 $L^1(\widetilde G)/L_{{\widetilde\Lambda}^\prime}^1(\widetilde G)$
  does not
 contain any copy of $\ell^1$, we have
 $$L(L^1(\widetilde G)/L_{{\widetilde\Lambda}^\prime}^1 (\widetilde
G),X)=
    K(L^1(\widetilde G)/L_{{\widetilde\Lambda}^\prime}^1 (\widetilde
G),X)$$ and since $
K(L^1(\widetilde G)/L_{{\widetilde\Lambda}^\prime}^1 (\widetilde
G),X)
= C_{\widetilde\Lambda}(\widetilde G,X)$, the proof is complete.
\end{pf}

For the next result, we need to recall some definitions.

\begin{df}
A bounded linear
operator $D: L^1[0,1] \longrightarrow X$ is called a Dunford-Pettis
operator if $D$ sends weakly compact sets into norm compact sets.
\end{df}
\begin{df}
An operator $T: L^1[0,1] \longrightarrow X$ is said to be nearly
representable if $T\circ D$ is (Bochner) representable for every
Dunford-Pettis operator $D: L^1[0,1] \longrightarrow L^1[0,1]$.
\end{df}

The following class of Banach spaces was introduced by Kaufman, Petrakis,
Riddle and Uhl in \cite{KPRU}:

\begin{df} A Banach space $X$ is said to have the Near Radon-Nikodym
Property (NearRNP) if every nearly representable operator from $L^1[0,1]$
into $X$ is representable.
\end{df}

\begin{thm} Let $X$ and $Y$ be Banach spaces such that:
\begin{itemize}
\item[(i)] $X$ is separable;
\item[(ii)] $X^*$ and $Y$ has the NearRNP;
\item[(iii)] $L(X,Y)=K(X,Y)$.
\end{itemize}
Then $K(X,Y)$ has the NearRNP.
\end{thm}

\begin{pf} Let $T: L^1[0,1] \longrightarrow K(X,Y)$ be a nearly
representable operator and let $F: \Sigma_{[0,1]} \longrightarrow K(X,Y)$
be the representing measure of the operator $T$. For each $x \in X$,
define $T_x: L^1[0,1] \longrightarrow Y$ as follows:
   $$T_x(f) = \langle T(f), x \rangle ,\quad \forall f \in L^1[0,1].$$
$T_x$ is clearly nearly representable and therefore representable (since
$Y$ has the NearRNP). Hence the measure $F(.)x$ which can be easely
checked to be the representing measure of $T_x$ has Bochner integrable
density. Similarly for each $y^* \in Y^*$, the operator $T^{y^*}:
L^1[0,1]
\longrightarrow X^*$ given by
$$T^{y^*}(f) = \langle (Tf)^*,y^* \rangle \quad \forall f \in L^1[0,1] $$
is nearly representable and therefore representable (since $X^*$ has the
NearRNP) and since the measure $F(.)^*y^*$ is the representing measure of
$T^{y^*}$, it has Bochner integrable density. Now we can procede as in
the proof of Theorem 2 to conclude that the measure $F$ has Bochner
integrable density which shows that the operator $T$ is representable.
\end{pf}

\noindent {\bf Remark:} Corollary 1,  Corollary 2 and Corollary 3
 are still valid if we replace the
$\Lambda$-RNP by  the NearRNP.

Let us now turn our attention to measures that can be represented by
Pettis-integrable functions. Recall that for a Banch space $X$, a
function
$f: G \longrightarrow X$ is Pettis-integrable if $f$ is weakly scalarly
measurable i.e. for each $x^* \in X^*$ the function $\langle f(.),x^*
\rangle$ is measurable
and for each $A \in \cal{B}(G)$, there exists $x_A \in X$ such that
 $$\langle x_A, x^* \rangle = \int_{A} \langle f(\omega),x^*\rangle \
d\lambda(\omega). $$
For more details about Pettis integral we refer to \cite{T4}.

\begin{df} Using the same notation as before, we say  that a Banach
space $X$ has type I-$\Lambda$-WeakRNP (resp. type II-$\Lambda$-WeakRNP)
if every measure $F$ in $M_{\Lambda}^\infty(G,X)$ (resp.
$M_{\Lambda,ac}(G,X)$) has Pettis-integrable density.
\end{df}

\begin{thm} Let $X$ and $Y$ be Banach spaces such that:
 \begin{itemize}
\item[(i)] $Y$ is a dual Banach space;
\item[(ii)] $X^*$ and $Y$ have the type I-$\Lambda$-WeakRNP (resp. type
II-$\Lambda$-WeakRNP);
\item[(iii)] $L(X,Y)=K(X,Y)$.
\end{itemize}
Then  $K(X,Y)$ has the type I-$\Lambda$-WeakRNP (resp. type
II-$\Lambda$-WeakRNP).
\end{thm}
\begin{pf}  Again we will present the type I case: Let $F: \cal{B}(G)
\longrightarrow K(X,Y)$ be a measure such that:
   \begin{itemize}
  \item[a)] $|F| \leq \lambda$;
  \item[b)] $\hat{F}(\gamma) =0 $ \ for $\gamma \notin \Lambda$.
  \end{itemize}
We will show that the measure $F$ has Pettis integrable density.
Notice first that the range of $F$ is separable. Using similar argument
as in section 1, we can assume without loss of generality that the
predual of $Y$ is separable.

Let $\rho$ be a lifting of $L^\infty(\lambda)$( see \cite{DIU} or
\cite{IT} for the
definition). Since $L(X,Y)=L(X,Z^*)$ is a dual space, we can apply
Theorem 4 of \cite{DIU}(P. 263) to get a
unique bounded function $T: G \longrightarrow L(X,Z^*)$ such that:
\begin{itemize}
 \item[(1)] $ \langle z, F(A)x \rangle = \dint_A \langle z, T(\omega)x
\rangle \ d\lambda(\omega)$ for $z \in Z$, $x \in X$ and $A \in
\cal{B}(G)$;
\item[(2)] $\rho(
\langle z, T(\omega)x \rangle ) =
\langle z, T(\omega)x \rangle $ for all $z \in Z$ and $x \in X$.
\end{itemize}
We claim that $\omega \longrightarrow T(\omega)$ is Pettis integrable.
We need several steps.

\begin{lem} For every $x \in X$ and $y^* \in Y^*$, the maps  $\omega
\longrightarrow T(\omega)x$  and $\omega \longrightarrow
T(\omega)^*y^*$  are Pettis-integrable and for  every $A \in \cal{B}(G)$
we have:
 \begin{itemize}
\item[(a)]  $F(A)x= \text{Pettis}-\dint_A T(\omega)\ d\lambda(\omega)$;
\item[(b)]  $F(A)^*y^*= \text{Pettis}-\dint_A T(\omega)^*y^* \
d\lambda(\omega)$.
\end{itemize}
\end{lem}

The proof can be done with essentially the same idea as in the proof
of Proposition 3, so we will omit the detail.

Let us now consider $\Delta=$ the unit ball of $Y^*$ with the
weak$^*$-topology. As in the proof of Theorem 3, we denote by $J$ the
natural isometry from $Y$ into $C(\Delta)$ and $J^\#$ the isometry from
$K(X,Y)$ into $K(X,C(\Delta))$.  Consider $\phi(\omega)=
J^\#(T(\omega))=J\circ T(\omega) \in K(X,C(\Delta))$. Since
$M(\Delta)=C(\Delta)^*$ has the metric approximation property, there
exists
a sequence of finite rank operators $\theta_n$ in
$L(C(\Delta),C(\Delta))$
such that $\sup\limits_{n \in \N}||\theta_n|| \leq 1$ and
$\theta_n$
converges to  the identity operator  uniformly on every compact subset of
$C(\Delta)$. Since the $\phi(\omega)$'s are compacts, we have
  $$\lim_{n\to \infty} ||\theta_n \circ \phi(\omega) -\phi(\omega)||=0$$
for every $\omega \in G$.

Let $\theta_n = \sum\limits_{k=1}^{k_n} \mu_{k,n} \otimes f_{n,k}$ where
$\mu_{k,n}$ and $f_{k,n}$ belong to $M(\Delta)$ and $C(\Delta)$
respectively. We have for every $\omega \in G$,
   $$\theta_n \circ \phi(\omega)x = \sum_{k=1}^{k_n} \langle
\phi(\omega)^* \mu_{k,n},x \rangle f_{k,n}$$
and if $I \in K(X,C(\Delta))^* = (X^* \hat{\otimes}_\epsilon
C(\Delta))^*= C(\Delta,X^*)^* = I(C(\Delta),X^{**})$; here
$I(C(\Delta),X^{**})$ denotes the  space of integral operators from
$C(\Delta)$ into $X^{**}$ (see \cite{DU} P. 232).
 We get
$$\langle \theta_n \circ \phi(\omega) , I \rangle =
 \sum_{k=1}^{k_n} \langle \phi(\omega)^*\mu_{k,n}, I( f_{k,n}) \rangle $$
which is measurable and therefore the map $\omega \longrightarrow
\theta_n \circ \phi(\omega)$ is weakly scalarly measurable and for every
$A \in \cal{B}(G)$ and $I \in K(X,C(\Delta))$, we have
      $$ \langle \theta_n \circ J \circ F(A),I \rangle
         =  \int_A \langle\theta_n \circ \phi(\omega), I \rangle \
d\lambda(\omega).$$
And now since $||\theta_n||\leq 1$, the set $\{\langle \theta_n \circ
\phi (.), I \rangle;\ n \in \N;\ ||I||\leq 1\}$ is uniformly
integrable so by taking the limit as
$n$ tends to
$\infty$, we get
(see
\cite{T4} Theorem 5-3-1) that
    $$J^\#(F(A))= \text{Pettis}-\int_A J^\#(T(\omega))\
d\lambda(\omega)$$
and since $J^\#$ is an isometry, the adjoint $(J^\#)^*$ is onto and
therefore the map $\omega \longrightarrow T(\omega)(G \longrightarrow
K(X,Y))$ is weakly scalarly measurable and it is now clear that $\omega
\longrightarrow T(\omega)$ is Pettis integrable with
$F(A)= \text{Pettis}-\dint_A T(\omega)\ d\lambda(\omega)$.
The proof is complete.
\end{pf}

\noindent{\bf Remark:}
 In \cite{EM3}, Emmanuele obtained the usual WeakRNP case of Theorem 4
but his method of proof is quite different and cannot be extended to the
general case of type I-$\Lambda$-WeakRNP or type II-$\Lambda$-WeakRNP.

Let us finish by asking the following question

\noindent
{\bf Question:}  If $X$ and $Y$ have type I-$\Lambda$-RNP (resp. type
II-$\Lambda$-RNP, resp. NearRNP, resp. WeakRNP) and
$L_{w^*}(X^*,Y)=K_{w^*}(X^*,Y)$; does $K_{w^*}(X^*,Y)$ have the same
property ? (here $L_{w^*}(X^*,Y)$(resp. $K_{w^*}(X^*,Y)$) denotes the
space of bounded (resp. compact bounded) operators from $X^*$ into $Y$
that are weak$^*$ to weak continuous).

The answer to this question is still unknown even for the usual RNP.

\noindent
{\bf Acknowlegements:}  The author would like to thank Professor Paula
Saab for several useful discussions conserning this paper.

\end{document}